\documentclass[runningheads, times]{llncs}
\usepackage{graphicx}
\usepackage{lineno,hyperref}
\usepackage{multirow}
\usepackage{rotating}
\usepackage{epstopdf}
\usepackage{graphicx}
\usepackage{subfigure}
\usepackage{threeparttable}
\usepackage{times}
\usepackage{booktabs} 

\begin{document}
\title{Adaptive RBF Interpolation for Estimating Missing Values in 
	Geographical Data
	\thanks{Supported by the 
		National Natural Science Foundation of China.}}
%
%
\author{Kaifeng Gao\inst{1} \and
Gang Mei\inst{1}* \and
Salvatore Cuomo \inst{2} \and
Francesco Piccialli \inst{2} \and
Nengxiong Xu \inst{1}
}
\authorrunning{K. Gao et al.}
%
\institute{School of Engineering and Technology, China University of Geosciences (Beijing), 100083, Beijing, China \and
	Department of Mathematics and Applications “R. Caccioppoli”, 
	University of Naples Federico II, Naples, Italy
	\\
	\email{gang.mei@cugb.edu.cn}\\
}
\maketitle              
\begin{abstract}
The quality of datasets is a critical issue in big data mining. More 
interesting things could be mined from datasets with higher quality. The 
existence of missing values in geographical data would worsen the quality of 
big datasets. To improve the data quality, the missing values are generally 
needed to be estimated using various machine learning algorithms or 
mathematical methods such as approximations and interpolations. In 
this paper, we propose an adaptive Radial Basis Function (RBF) interpolation 
algorithm for estimating missing values in geographical data. In the 
proposed method, the samples with known values are considered as the data 
points, while the samples with missing values are considered as the 
interpolated points. For each interpolated point, first, a local set of data 
points are adaptively determined. Then, the missing value of the 
interpolated point is imputed via interpolating using the RBF interpolation 
based on the local set of data points. Moreover, the shape factors of the 
RBF are also adaptively determined by considering the distribution of the 
local set of data points. To evaluate the performance of the proposed 
method, we compare our method with the commonly used $k$ Nearest Neighbors ($k$NN) interpolation and 
Adaptive Inverse Distance Weighted (AIDW) methods, and conduct three groups of benchmark experiments. Experimental 
results indicate that the proposed method outperforms the $k$NN interpolation 
and AIDW in terms of accuracy, but worse than the $k$NN interpolation and AIDW 
in terms of efficiency.

\keywords{Data mining \and Data quality \and Data imputation \and RBF 
	interpolation \and $k$NN}
\end{abstract}
\section{Introduction}
Datasets are the key elements in big data mining, and the quality of 
datasets has an important impact on the results of big data analysis. For a 
higher quality dataset, some hidden rules can often be mined from it, and 
through these rules we can find some interesting things. At present, big 
data mining technology is widely used in various fields, such as geographic 
analysis \cite{17,13}, financial analysis, smart city and biotechnology. It 
usually needs a better dataset to support the research, but in fact there is 
always noise data or missing value data in the datasets \cite{15,11,14}. In order to improve 
data quality, various machine learning algorithms \cite{16,12} are often required to 
estimate the missing value and clean noise data.

RBF interpolation algorithm is a popular method for estimating missing 
values \cite{2,3,1}.In large-scale computing, the cost can be minimized by using adaptive scheduling method \cite{18}.
 RBF is a distance-based function, which is meshless and 
dimensionless, thus it is inherently suitable for processing 
multidimensional scattered data. Many scholars have done a lot of work on 
RBF research. Skala \cite{9} used CSRBF to analyze big datasets, Cuomo et al. 
\cite{6,7} studied the reconstruction of implicit curves and surfaces by RBF 
interpolation. Kedward et al. \cite{8} used multiscale RBF interpolation to study 
mesh deformation. In RBF, the shape factor is an important reason affecting 
the accuracy of interpolation. Some empirical formulas for optimum shape 
factor have been proposed by scholars.

In this paper, our objective is to estimate missing values in geographical 
data.We proposed an adaptive RBF interpolation algorithm, which adaptively 
determines the shape factor by the density of the local dataset. To evaluate 
the performance of adaptive RBF interpolation algorithm in estimating 
missing values, we used three datasets for verification experiments, and 
compared the accuracy and efficiency of adaptive RBF interpolation with that 
of $k$NN interpolation and AIDW.

The rest of the paper is organized as follows. Section 2 mainly introduces 
the implementation process of the adaptive RBF interpolation algorithm, and 
briefly introduces the method to evaluate the performance of adaptive RBF 
interpolation. Section 3 introduces the experimental materials. Section 4 
presents the estimated results of missing values, then discusses the 
experimental results. Section 5 draws some conclusions.

\section{Methods}
In this paper, our objective is to develop an adaptive RBF interpolation 
algorithm to estimate missing values in geospatial big data, and compare the 
results with that of $k$NN and AIDW. In this section, we firstly introduce the 
specific implementation process of the adaptive RBF interpolation algorithm, 
then briefly introduces the method to evaluate the performance of adaptive 
RBF interpolation.

\subsection{Adaptive RBF Interpolation Algorithm}

The basic ideas behind the RBF interpolation are as follows. Constructing a 
high-dimensional function $f\left( x \right),x\in R^n$£¬suppose there is a 
set of discrete points $x_i \in R^n,i=1,2,\cdots N_{£¬}$with associated 
data values $f\left( {x_i } \right)\in R,i=1,2,\cdots N$ Thus, the function 
$f\left( x \right)$ can be expressed as a linear combination of RBF in the 
form (Eq. (\ref{eq1})):
\begin{equation}
\label{eq1}
f\left( x \right)=\sum\nolimits_{j=1}^N {a_j \phi \left( {\left\| {x-x_j } 
		\right\|_2 } \right)} 
\end{equation}
where $N$ is the number of interpolation points $\left\{ {a_j } \right\}$ is 
the undetermined coefficient, the function $\phi $ is a type of RBF 

The kernel function selected in this paper is Multi-quadric RBF(MQ-RBF), 
which is formed as (Eq. (\ref{eq2})):
\begin{equation}
\label{eq2}
\phi \left( r \right)=\sqrt {\left( {r^2+c^2} \right)} 
\end{equation}
where $r$ is the distance between the interpolated point and the data point, 
$c$ is the shape factor. Submit the data points $\left( {x_i ,y_i } \right)$ 
into Eq. (\ref{eq1}), then the interpolation conditions become (Eq. (\ref{eq3})):
\begin{equation}
\label{eq3}
y_i =f\left( {x_i } \right)=\sum\nolimits_{j=1}^N {a_j \phi \left( {\left\| 
		{x_i -x_j } \right\|_2 } \right)} ,i=1,2,\cdots N
\end{equation}

When using the RBF interpolation algorithm in a big dataset, it is not 
practical to calculate an interpolated point with all data points. 
Obviously, the closer the data point is to the interpolated point, the 
greater the influence on the interpolation result and the data point far 
from the interpolated point to a certain distance, its impact on the 
interpolated point is almost negligible. Therefore, we calculate the 
distances from an interpolated point to all data points, and select 20 
points with the smallest distances as a local dataset for the interpolated 
point.

In Eq. (\ref{eq2}), the value of the shape factor $c$ in MQ-RBF has a significant 
influence on the calculation result of interpolation. We consult the method 
proposed by Lu and Wang \cite{4,5}, adaptively determining the value $c$ of the 
interpolated points by the density of the local dataset. The expected 
density $D_{\exp } $ is calculated by the function (Eq. (\ref{eq4})):
\begin{equation}
\label{eq4}
D_{\exp } =\frac{N_{dp} }{\left( {X_{\max } -X_{\min } } \right)\left( 
	{Y_{\max } -Y_{\min } } \right)}\quad
\end{equation}
where $N_{dp} $ is the number of data points in the dataset, $X_{\max } $ is 
the maximum value of $x_i $ for the data points in the dataset, $X_{\min } $ 
is the minimum value of $x_i $ in dataset $Y_{\max } $ is the maximum value 
of $y_i $ in dataset, $Y_{\min } $ is the minimum value of $y_i $ in dataset

And the local density $D_{loc} $ is calculated by (Eq. (\ref{eq5})):
\begin{equation}
\label{eq5}
D_{loc} =\frac{N_{loc} }{\left( {x_{\max } -x_{\min } } \right)\left( 
	{y_{\max } -y_{\min } } \right)}
\end{equation}
where $N_{loc} $ is the number of data points in the local dataset, in this 
paper, we set $N_{loc} $ as 20. $x_{\max } $ is the maximum value of $x_i $ 
for the data points in local dataset, $x_{\min } $is the minimum value of 
$x_i $ in local dataset $y_{\max } $ is the maximum value of $y_i $ in local 
dataset, $y_{\min } $ is the minimum value of $y_i $ in local dataset

With both the local density and the expected density, the local density 
statistic $D$ can be expressed as (Eq. (\ref{eq6})):
\begin{equation}
\label{eq6}
D\left( {s_0 } \right)=\frac{D_{loc} }{D_{\exp } }
\end{equation}
where $s_0 $ is the location of an interpolated point Then normalize the 
$D\left( {s_0 } \right)$ measure to $\mu _D $ by a fuzzy membership 
function (Eq. (\ref{eq7})):
\begin{equation}
\label{eq7}
_{\mu _D =\left\{ {{\begin{array}{*{20}c}
			0 & {D\left( {s_0 } \right)\le 0} \\
			{0.5-0.5\cos \left[ {\frac{\pi }{2}D\left( {s_0 } \right)} \right]} & {0\le 
				D\left( {s_0 } \right)\le 2} \\
			1 & {D\left( {s_0 } \right)\ge 2} \\
			\end{array} }} \right.} 
\end{equation}
Then determine the shape factor $c$ by a triangular membership function. See 
Eq (\ref{eq8}).
\begin{equation}
\label{eq8}
_{c = \left\{ {{\begin{array}{*{20}c}
			{c_1 } & {0.0\le \mu _D \le 0.1} \\
			{c_1 \left[ {1-5\left( {\mu _D -0.1} \right)} \right]+5c_2 \left( {\mu _D 
					-0.1} \right)} & {0.1\le \mu _D \le 0.3} \\
			{5c_3 \left( {\mu _D -0.3} \right)+c_2 \left[ {1-5\left( {\mu _D -0.3} 
					\right)} \right]} & {0.3\le \mu _D \le 0.5} \\
			{c_3 \left[ {1-5\left( {\mu _D -0.5} \right)} \right]+5c_4 \left( {\mu _D 
					-0.5} \right)} & {0.5\le \mu _D \le 0.7} \\
			{5c_5 \left( {\mu _D -0.7} \right)+c_4 \left[ {1-5\left( {\mu _D -0.7} 
					\right)} \right]} & {0.7\le \mu _D \le 0.9} \\
			{c_5 } & {0.9\le \mu _D \le 1.0} \\
			\end{array} }} \right.}\quad
\end{equation}
where $c_1 ,c_2 ,c_3 ,c_4 ,c_5 $ are five levels of shape factor.

\begin{figure}[htbp]
	\centering
	\includegraphics[width=0.9\textwidth]{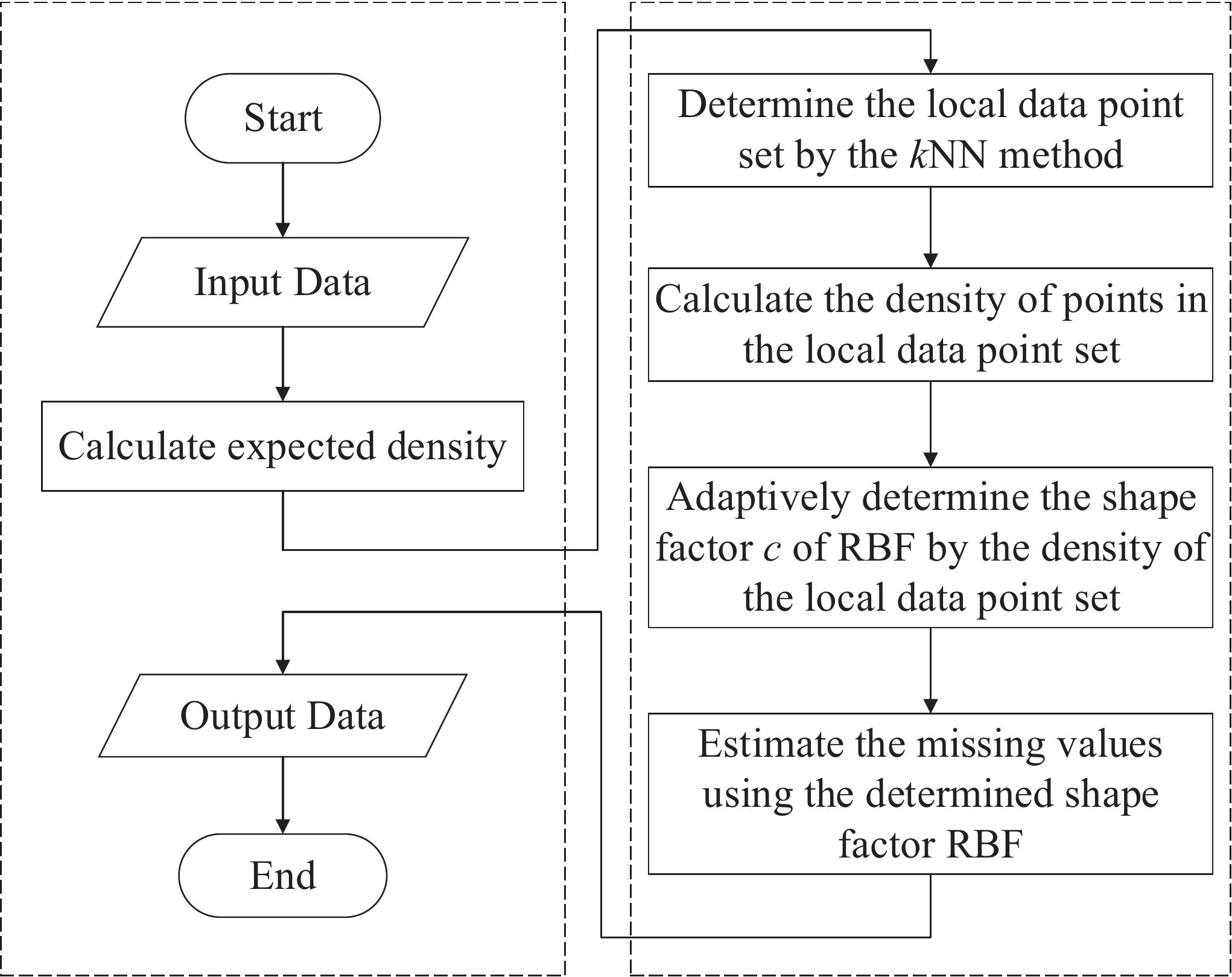}
	\caption{Flowchart of the adaptive RBF interpolation algorithm } \label{fig1}
\end{figure}

After determining the shape factor $c$, the next steps are the same as the 
general RBF calculation method. The specific process of the adaptive RBF 
interpolation algorithm is illustrated in Fig.\ref{fig1}. 

\subsection{Evaluating the Performance of Adaptive RBF Interpolation}

In order to evaluate the computational accuracy of the adaptive RBF 
interpolation algorithm, we use the metric, Root Mean Square Error (RMSE) to 
measure the accuracy. The RMSE evaluates the error accuracy by comparing the 
deviation between the estimated value and the true value. Then, we compare 
the accuracy and efficiency of adaptive RBF estimator with the results of 
$k$NN and AIDW estimators.

\section{Experimental Design}


To evaluate the performance of the presented adaptive RBF interpolation 
algorithm we use three datasets to test it. The details of the experimental 
environment are listed in Table \ref{tab1}.

\begin{figure}[!h]
	\centering
	\subfigure[DEM map of Beijing City]{
		\label{fig2}       
		\includegraphics[width=0.4\textwidth]{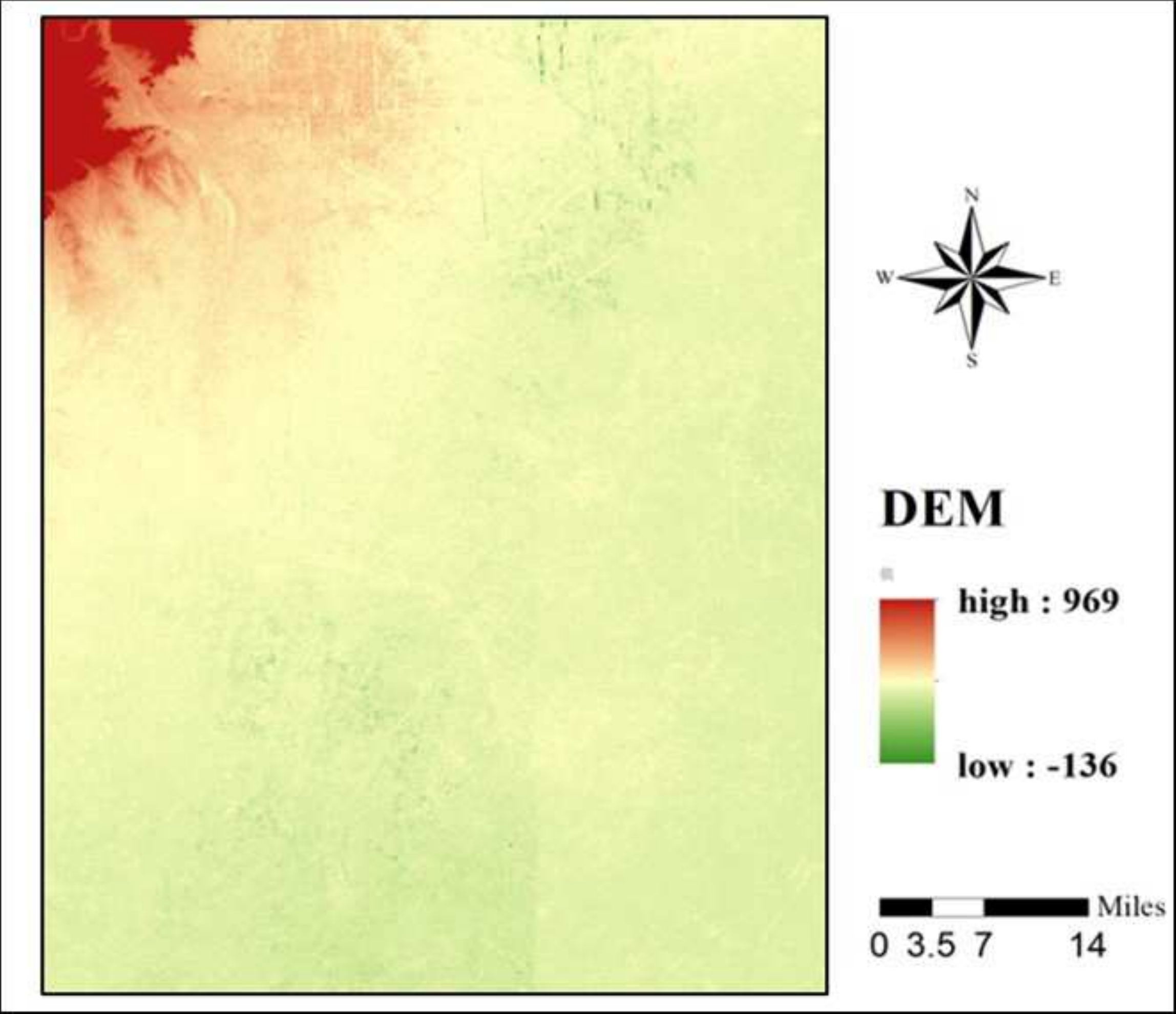}
	}
	\hspace{1em}
	\subfigure[DEM map of Chongqing City]{
		\label{fig3}       
		\includegraphics[width=0.415\textwidth]{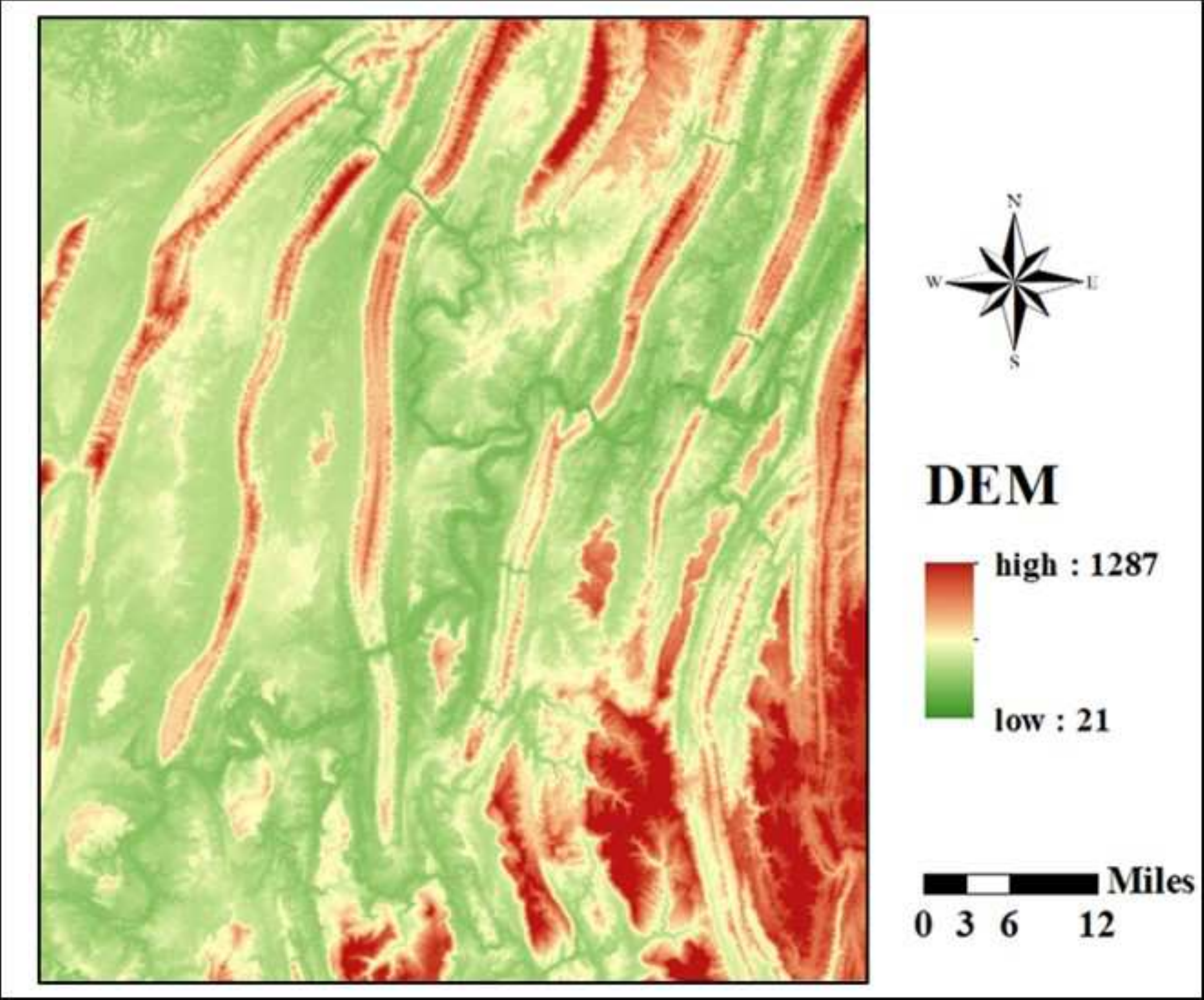}
	}
	\hspace{1em}
	\subfigure[DEM map of Longyan City]{
		\label{fig4}       
	    \includegraphics[width=0.415\textwidth]{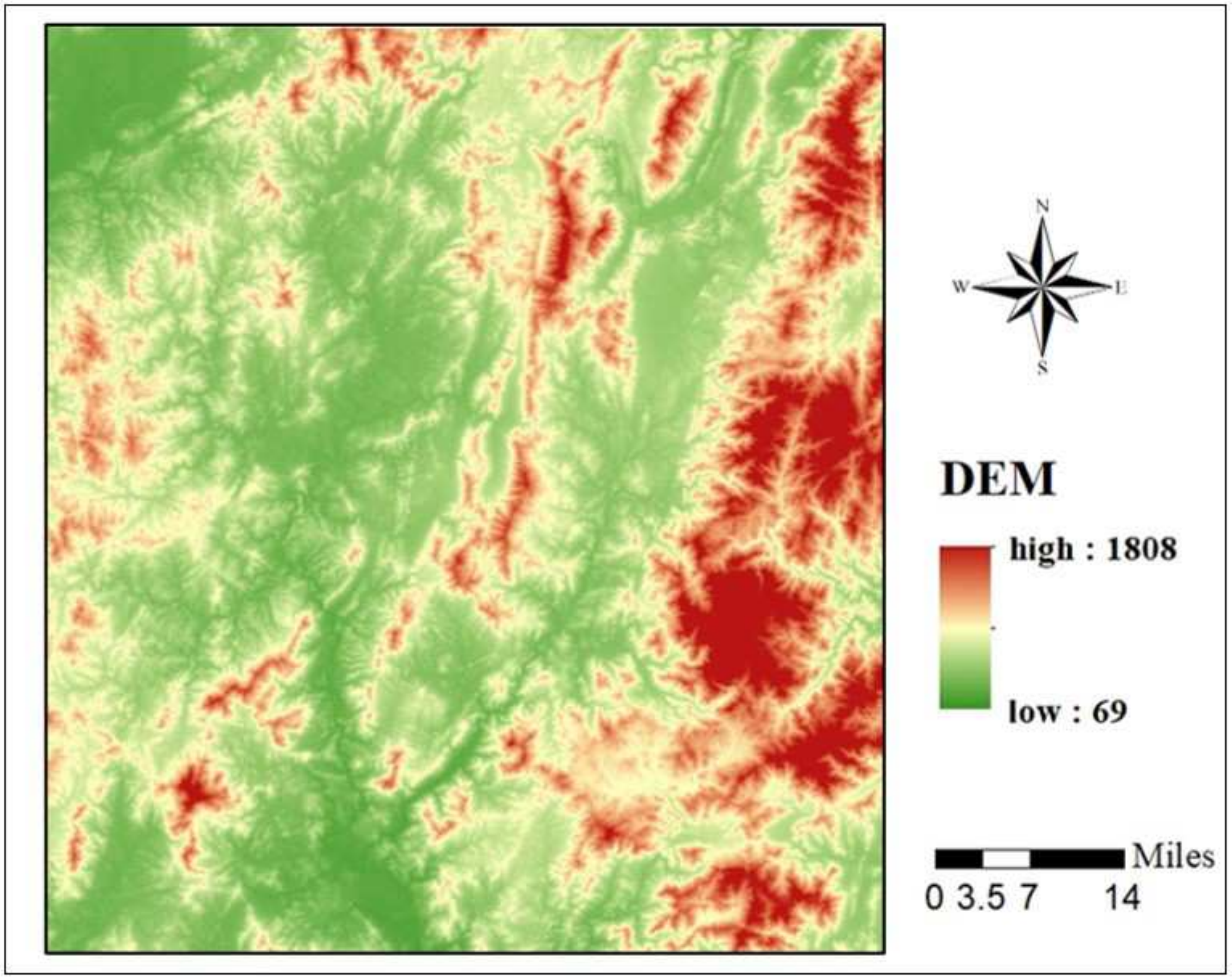}
	}
	\caption{The DEM maps of three cities for the experimental tests.}
\label{fig:5}       
	\label{figs234}       
\end{figure}

In our experiments, we use three datasets from three cities' digital 
elevation model (DEM) (Fig.\ref{fig2} to Fig.\ref{fig4}), the range of three DEM maps is 
the same We randomly select 10{\%} observed samples from each dataset as the 
samples with missing values, and the rest as the samples with known values. 
It should be noted that the samples with missing values have really 
elevation values in fact, but for testing, we assume the elevations are 
missing. Basic information of the datasets is listed in Table \ref{tab2}.

\begin{table}[htbp]
	\caption{Details of experimental environment}
	\begin{center}
		\begin{tabular}{|l|l|}
			\hline
			\textbf{Specification}& 
			\textbf{Details} \\
			\hline
			\textbf{OS}& 
			Windows 7. Professional \\
			\hline
			\textbf{CPU}& 
			Intel (R) i5-4210U \\
			\hline
			\textbf{CPU Frequency}& 
			1.70 GHz \\
			\hline
			\textbf{CPU RAM}& 
			8 GB \\
			\hline
			\textbf{CPU Core}& 
			4 \\
			\hline
		\end{tabular}
		\label{tab1}
	\end{center}
\end{table}

\begin{table}[htbp]
	\caption{Details of the testing data}
	\begin{center}
		\begin{tabular}{|l|l|l|l|}
			\hline
			\textbf{Dataset}& 
			\textbf{Number of known values}& 
			\textbf{Number of missing values}& 
			\textbf{Illustration} \\
			\hline
			\textbf{Beijing}& 
			1,111,369& 
			123,592& 
			Fig.\ref{fig2} \\
			\hline
			\textbf{Chongqing}& 
			1,074,379& 
			97,525& 
			Fig.\ref{fig3} \\
			\hline
			\textbf{Longyan}& 
			1,040,670& 
			119,050& 
			Fig.\ref{fig4} \\
			\hline
		\end{tabular}
		\label{tab2}
	\end{center}
\end{table}

\section{Results and~Discussion}

We compare the accuracy and efficiency of adaptive RBF estimator with that 
of $k$NN and AIDW estimators.

\begin{figure}[htbp]
	\centering
	\includegraphics[width=0.75\textwidth]{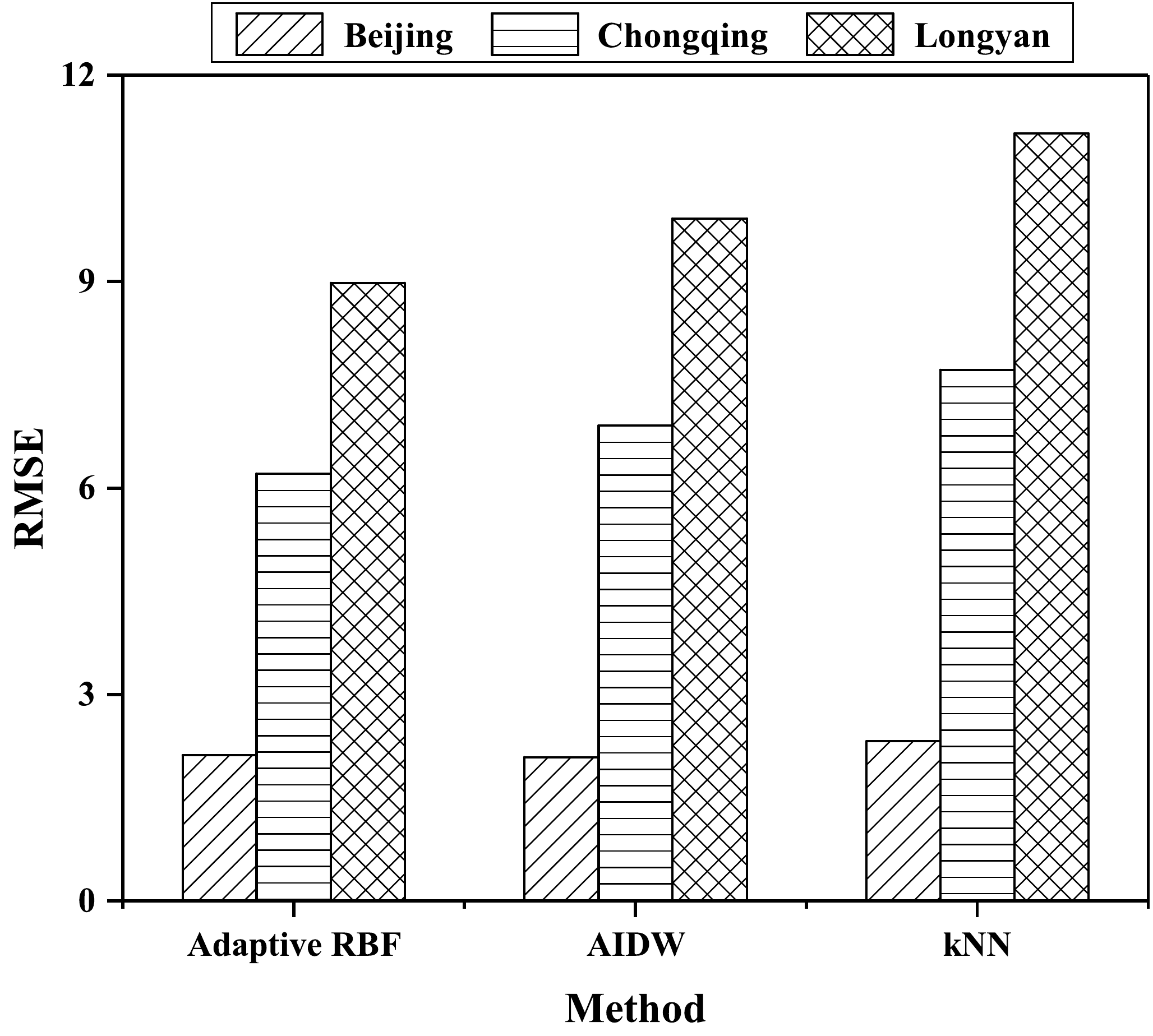}
	\caption{Comparisons of the computational accuracy} \label{fig5}
\end{figure}

\begin{figure}[htbp]
	\centering
	\includegraphics[width=0.75\textwidth]{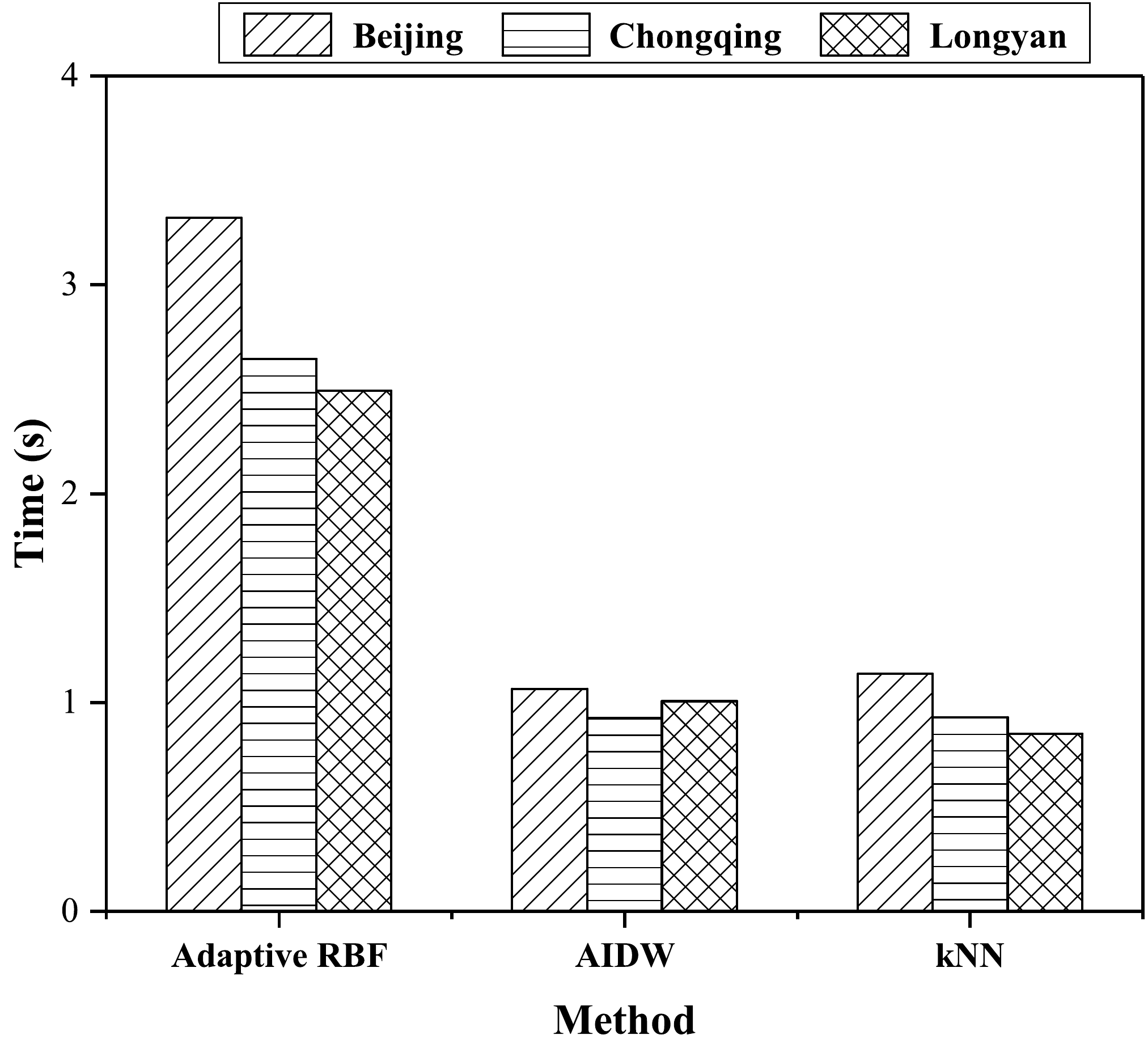}
	\caption{Comparisons of the computational efficiency} \label{fig6}
\end{figure}

In Fig.\ref{fig5} and Fig.\ref{fig6}, we find that the accuracy of the adaptive RBF estimator is 
the best performing, and the$ k$NN estimator with the lowest accuracy. With the 
number of known data points in the datasets decreases, the accuracy of three 
estimators decreases significantly. Moreover, the computational efficiency 
of adaptive RBF estimator is worse than that of $k$NN estimator and AIDW 
estimator, among them,$ k$NN has the best computational efficiency. With the 
increase of data quantity, the disadvantage of the computational efficiency 
of $k$NN estimator becomes more and more obvious.


The data points selected from DEM are evenly distributed, and the shape 
factor $c$ of the adaptive RBF interpolation algorithm is adapted according to 
the density of the points in the local dataset therefore, when the missing 
data is estimated in a dataset with a more uniform data point, the 
advantages of the adaptive RBF interpolation algorithm may not be realized. 
We need to do further research in datasets with uneven datasets 

\section{Conclusions}

In this paper, we specifically proposed an adaptive RBF interpolation algorithm for 
estimating the missing values in geographical data. We performed three groups of 
experimental tests to evaluate the computational accuracy and efficiency of the 
proposed adaptive RBF interpolation by comparing with the  $k$NN interpolation and AIDW method. We found that the accuracy of the 
adaptive RBF interpolation performs better than $k$NN interpolation and AIDW in 
regularly distributed datasets. But the efficiency of adaptive RBF 
interpolation is worse than the $k$NN interpolation and AIDW.

\subsubsection*{Acknowledgments}
This research was jointly supported by the 
National Natural Science Foundation of China (11602235), and the Fundamental 
Research Funds for China Central Universities (2652018091, 2652018107, and 2652018109).

%
%
\bibliographystyle{splncs04}
\bibliography{References}

\end{document}